\documentstyle[latexsym,amsfonts]{elsart}

\input{amssym.def}
\input{amssym.tex}

\newcommand{\N}{{\mathbb N}}

\newcommand{\R}{{\mathbb R}}

\newcommand{\Sp}{{\mathbb S}}
\newcommand{\eqref}[1]{(\ref{#1})}
\newcommand{\supdt}{\sup_{0\le\tau\le t}\|d(\tau)\|}  
\newcommand{\supdT}{\sup_{0\le\tau\le T}\|d(\tau)\|}  

\def\KK{{\cal K}}

\def\LL{{\cal L}}
\def\DD{{\cal D}}

\def\frall{{\rm\enspace for\enspace all\enspace}}
\newcommand{\norm}[1]{\left\Vert#1\right\Vert}

\newcounter{enumctr}

\begin{document}
\begin{frontmatter}
  \title{Asymptotic stability equals exponential stability, and ISS
    equals finite energy gain---if you twist your eyes}
   \thanks[lars]{This paper has been written while the first
    author was visiting the Dipartimento di Matematica, Universit\'a di
Roma 
       ``La Sapienza'', Italy, supported by 
       DFG-Grant GR1569/2-1.}
    \thanks[eds]{Supported in part by US Air Force Grant
F49620-98-1-0242}
\author{Lars Gr\"une  \thanksref{lars}}
\address{Fachbereich Mathematik,  J.W. Goethe-Universit\"at,
Postfach 11 19 32, D-60054 Frankfurt a.M., Germany,
E-Mail:~gruene@math.uni-frankfurt.de}
\author{Eduardo D.~Sontag  \thanksref{eds}}
\address{Department of Mathematics, Rutgers University, New Brunswick,
NJ 08903, USA, E-Mail:~sontag@control.rutgers.edu}
\author{Fabian R.~Wirth}
\address{Zentrum f\"ur Technomathematik,  Universit\"at Bremen,
  D-28344 Bremen, Germany, E-Mail:~fabian@math.uni-bremen.de}

  \begin{abstract}
  In this paper we show that uniformly global asymptotic stability for a 
  family of ordinary differential equations is equivalent to uniformly
  global exponential stability under a suitable nonlinear change of
variables. 
  The same is shown for input-to-state stability and 
  input-to-state exponential stability, and for 
  input-to-state exponential stability and a nonlinear $H_\infty$
estimate.
  \end{abstract}

\end{frontmatter}
\section{Introduction}

Lyapunov's notion of (global) asymptotic stability of an equilibrium
is a key concept in the qualitative theory of differential equations
and nonlinear control.  In general, a far stronger property is that of
{\em exponential\/} stability, which requires decay estimates of the
type ``$\norm{x(t)}\leq ce^{-\lambda t}\norm{x(0)}$.''  (See for
instance \cite{Khalil95} for detailed discussions of the comparative
roles of asymptotic and exponential stability in control theory.)  In
this paper, we show that, for differential equations evolving in
finite-dimensional Euclidean spaces $\R^n$ (at least in spaces of
dimensions $\not= 4,5$) the two notions are one and the same under
coordinate changes.

Of course, one must define ``coordinate change'' with care, since
under diffeomorphisms the character of the linearization at the
equilibrium (which we take to be the origin) is invariant.  However,
if, in the spirit of both structural stability and the classical
Hartman-Grobman Theorem (which, cf.~\cite{Perko}, gives in essence a
{\em local\/} version of our result in the special hyperbolic case),
we relax the requirement that the change of variables be {\em smooth
  at the origin\/}, then all obstructions disappear.  Thus, we ask
that transformations be infinitely differentiable except possibly at
the origin, where they are just continuously differentiable. Their
respective inverses are continuous globally, and infinitely
differentiable away {}from the origin.

Closely related to our work is the fact that all asymptotically stable
{\em linear\/} systems are equivalent (in the sense just discussed) to
$\dot{x}=-x$; see e.g.~\cite{Arno92}.  The basic idea of the proof
in~\cite{Arno92} is based upon projections on the level sets of
Lyapunov functions, which in the linear case of course be taken to be
quadratic (and hence have ellipsoids as level sets).  It is natural to
use these ideas also in the general nonlinear case, and Wilson's
paper~\cite{Wils67}, often cited in control theory, remarked that
level sets of Lyapunov functions are always homotopically equivalent
to spheres.  Indeed, it is possible to obtain, in great generality, a
change of coordinates rendering the system in normal form $\dot x =
-x$ (and hence exponentially stable), and several partial versions of
this fact have appeared in the literature, especially in the context
of generalized notions of homogeneity for nonlinear systems; see for
instance~\cite{Cole65,PomePral88,Kawski,Rosier-thesis,Pral97}.

It is perhaps surprising that, at least for unperturbed systems, the
full result seems not to have been observed before, as the proof is a
fairly easy application of results from differential topology.  (Those
results are nontrivial, and are related to the generalized Poincar\'e
conjecture and cobordism theory; in fact, the reason that we only make
an assertion for $\not= 4,5$ is closely related to the fact that the
original Poincar\'e conjecture is still open.)

Note, however, that it has been common practice in the papers treating
the nonlinear case to use the flow generated by the original system to
define an equivalence transformation, thereby reducing the regularity
of the transformation to that of the system.  Here we use the flow
generated by the (normalized) Lyapunov function itself, which yields
more regular transformations.  In addition, and most importantly, our
poof also allows for the treatment of perturbed systems (for which the
reduction to $\dot x = -x$ makes no sense).

Lyapunov's notion is the appropriate generalization of exponential
stability to nonlinear differential equations.  For systems with
inputs, the notion of {\em input to state stability (ISS)\/}
introduced in~\cite{coprime} and developed further in
\cite{Teel,issmtns,Isidori,Jiang-Teel-Praly,Krstic-Deng,KKKbook,%
  PraW96,Sepulchre,SonW95,Tsinias} and other references, has been
proposed as a nonlinear generalization of the requirement of finite
${\cal L}^2$ gain or, as often also termed because of the spectral
characterizations valid for linear systems, ``finite nonlinear
$H^\infty $ gain'' (for which see e.g.\ 
\cite{Basar,Helton-James,Isidori-Astolfi,vdS}).  We also show in this
paper that under coordinate changes (now in both state and input
variables), the two properties (ISS and finite $H^\infty $ gain)
coincide (again, assuming dimension $\not= 4,5$).

We do not wish to speculate about the implications of the material
presented here.  Obviously, there are no ``practical'' consequences,
since finding a transformation into an exponentially stable system is
no easier than establishing stability (via a Lyapunov function).
Perhaps these remarks will be of some use in the further theoretical
development of ISS and other stability questions.  In any case, they
serve to further justify the naturality of Lyapunov's ideas and of
concepts derived from his work.

\section{Setup}

We consider the family of differential equations
\begin{equation}
\label{system}
  \dot{x}(t) = f(x(t),d(t))
\end{equation}
where $f: \R^n \times D \to \R^n $ is continuous and for $x\neq 0$
locally Lipschitz continuous in $x$, where the local Lipschitz
constants can be chosen uniformly in $d\in D\subseteq \R^m$.  Let
${\cal D}$ denote the set of measurable, locally essentially bounded
functions {}from $\R$ to $D$.  For any $x_0\in\R^n$ and any
$d(\cdot)\in\DD$,
there exists at least one maximal solution of~\eqref{system} for $t\geq0$,
with $x(0)=x_0$.
By abuse of notation, we denote any such solution, even if not unique,
as $\phi(t,x_0,d(\cdot))$, $t\in I(x,d(\cdot))$,
where $I(x,d(\cdot))$ is its existence interval.
Throughout the paper, $\|\cdot\|$ denotes the usual
Euclidean norm, and ``smooth'' means $C^\infty$. For a differentiable
function
$V:\R^n\to \R$ the expression $L_{f_d}V(x)$ denotes the directional
derivative
$DV(x)f(x,d)$.

The general framework afforded by the model~(\ref{system}) allows us
to treat simultaneously classical differential equations (the case
when $D=\{0\}$) and more generally robust stability of differential
equations subject to perturbations (when functions in ${\cal D}$ are
seen as disturbances which do not change the equilibrium, as in
parameter uncertainty), as well as systems with inputs in which
elements of ${\cal D}$ are seen as exogenous tracking or regulation
signals, or as actuator errors (in which case, the continuity
properties of $(x,d)\mapsto\phi(\cdot,x,d)$ are of interest).  In
light of these applications, we now describe the appropriate stability
concepts.

For the first, assume that $D$ is compact and that $f(0,d) =0$ for all
$d\in D$.  Then we say that the zero state is uniformly globally
asymptotically stable (UGAS) if there exists a class $\KK\LL$ function
$\beta$ such that, for each $d(\cdot)\in\DD$,
every maximal solution is defined for all $t\geq0$ and
\begin{equation} 
\|\phi(t,x,d(\cdot))\| \le \beta(\|x\|,t)
\label{gas} \end{equation}
for all $t\ge 0$.  As usual, we call a function
$\alpha:[0,\infty)\to [0,\infty)$ of class $\KK$, if it satisfies
$\alpha(0)=0$ and is continuous and strictly increasing (and class
$\KK_\infty$ if it is unbounded), and we call a continuous function
$\beta:[0,\infty)^2\to [0,\infty)$ of class $\KK\LL$, if it is
decreasing to
zero in the second and of class $\KK$ in the first argument.  (It is an
easy
exercise, cf.\ e.g.~\cite{LiSW96}, to verify that this definition is
equivalent to the requirements of uniform stability and uniform
attraction
stated in ``$\varepsilon-\delta$'' terms.) Note that while our general
assumptions on the right hand side $f$ do not guarantee uniqueness of
solutions through zero, the added assumption of asymptotic stability
implies
that $\phi(t,0,d)\equiv 0$ is the unique solution with initial condition
$x=0$, for all $d\in \DD$.
As a consequence, since away from zero we have a local Lipschitz condition,
solutions are unique for each given initial state and $d\in \DD$.

If the origin is no common fixed point for all values $d\in D$ then
\eqref{gas} is impossible. In this case, however, still a useful
notion of stability is possible. We call the system \eqref{system}
(globally) input-to-state stable (ISS), if there exists a class
$\KK\LL$ function $\beta$ and a class $\KK_\infty$ function $\alpha$
such that all solutions of \eqref{system} satisfy
\begin{equation} 
\|\phi(t,x,d(\cdot))\| \le \beta(\|x\|,t) + \alpha(\supdt) 
\label{iss} \end{equation}
for all $d(\cdot)\in\DD$ and all $t\ge 0$. Formulation \eqref{iss} is
the most frequently used characterization of the ISS property. Note
that with $\tilde{\beta}=2\beta$ and $\tilde{\alpha}=2\alpha$
inequality \eqref{iss} immediately implies
\[ \|\phi(t,x,d(\cdot))\| \le 
\max\left\{ \tilde{\beta}(\|x\|,t),\, \tilde{\alpha}(\supdt)\right\}\,
,\] 
hence this ``$\max$'' formulation can be used as an equivalent 
characterization.

Two apparently stronger formulations of these properties are obtained if
we
replace $\beta(\|x\|,t)$ by $ce^{-\lambda t}\|x\|$, more precisely we
call the
zero position of \eqref{system} uniformly globally {\em exponentially}
stable
(UGES), if there exist constants $c\ge 1, \lambda>0$ such that
\begin{equation} 
\|\phi(t,x,d(\cdot))\| \le ce^{-\lambda t}\|x\|
\label{ges} \end{equation}
holds for all $d(\cdot)\in\DD$ and all $t\ge 0$, and we call the system
input-to-state {\em exponentially} stable (ISES), if there exist a
class 
$\KK_\infty$ function $\alpha$ and constants $c\ge 1, \lambda>0$ such
that
\begin{equation} 
\|\phi(t,x,d(\cdot))\| \le \max\left\{ ce^{-\lambda t}\|x\|, 
\, \alpha(\supdt) \right\}
\label{ises} \end{equation}
for all $d(\cdot)\in\DD$ and all $t\ge 0$. (As usual, these definitions
use
appropriate constants $c, \lambda>0$.  In this paper, however, we will
see
that we can always work with ``normalized'' versions choosing $c=1,
\lambda=1$. For the (ISES) property we use the ``$\max$'' formulation
because
it allows a further implication as stated in Theorem 5, below. Observe
that
\eqref{ises} implies \eqref{iss} with $\beta(\|x\|,t) = ce^{-\lambda
  t}\|x\|$.)

Extending the concepts in \cite[p. 207]{Arno92} to our nonlinear setting,
we will call a homeomorphism
\[
T:\R^n \to\R^n
\]
a {\em change of variables\/}
if $T(0) = 0$, $T$ is $C^1$ on $\R^n$,
and $T$ is diffeomorphism on $\R^n\setminus \{0\}$
(i.e., the restrictions of $T$ and of $T^{-1}$ to $\R^n \setminus \{0\}$ are
both smooth).
Given a change of variables $T$ and a system~\eqref{system}, we may consider
the {\em transformed system\/}
\begin{equation}
  \label{transform}
  \dot{y}(t) = \tilde{f}(y(t),d(t)) \,,
\end{equation}
where, by definition,
\[
\tilde{f}(y,d) = DT(T^{-1}(y)) f(T^{-1}(y),d) \,.
\]
In other words, system~\eqref{transform} is obtained from the original system
by means of the change of variables $y=T(x)$.
Observe that the new system again satisfies the general requirements:
$\tilde{f}(y,d)$ is continuous, and it is locally Lipschitz on $x$
for $x\not=0$, uniformly on $d$.

It is our aim to show that for dimensions $n\neq 4,5$ the following
assertions are true.
Given a system of the form \eqref{system} satisfying \eqref{gas} or
\eqref{iss}, respectively, there exists a transformed system that
satisfies \eqref{ges} or \eqref{ises}, respectively.
In this sense, global asymptotic stability
is equivalent to global exponential stability under nonlinear changes of
coordinates. Furthermore, one may obtain transformed systems where the
constants defining the exponential stability property
can be chosen to be the special values $c=\lambda=1$.
 
Furthermore we show that if system \eqref{system} is ISES \eqref{ises}
with $c=\lambda=1$ 
then there exists a homeomorphism $R:\R^m \to \R^m$ on the input space with
$R(0) =
0$ that is a diffeomorphism on $\R^m \setminus \{0\}$ such that the
transformed system with $v=R(d)$
\begin{equation} \label{dtransform}
  \dot{y}(t) = \bar{f}(x(t),v(t)),\quad 
  \bar{f}(x,v) = f(x,R^{-1}(v))
\end{equation}
satisfies the following ``$L_2$ to $L_2$'' nonlinear $H_\infty$
estimate:
\begin{equation} \label{hinf}
\int_0^t \|\phi(s,x,v(\cdot))\|^2 ds \le \|x\|^2 + 
\int_0^t \|v(s)\|^2ds.
\end{equation}  
Since \eqref{hinf} in turn implies ISS 
(by \cite[Theorem 1]{Sont98}), we obtain equivalence between ISS 
and the nonlinear $H_\infty$ estimate \eqref{hinf} up to nonlinear
changes 
of coordinates.

\section{Construction of the coordinate transformation}
\label{Stab}

The main tool for our construction of $T$ is the use of an appropriate
Lyapunov function $V$. In fact, we can obtain $T$ for a whole class of
functions as stated in the following proposition. Recall that a
function $V:\R^n\to\R$ is called {\em positive definite} if $V(0)=0$
and $V(x)>0$ for all $x\ne 0$, and {\em proper} if the set $\{x\,|\,
V(x)\le \beta\}$ is bounded for each $\beta\ge 0$.

The next result says in particular that any such function may look
like $\|x\|^2$ under a coordinate change. This implies in particular
that the level sets under coordinate change are spheres. It may
therefore not come as a surprise that a basic ingredient of the proof
is related to the question of whether level sets of Lyapunov functions
in $\R^n$ are diffeomorphic to the sphere $S^{n-1}$. This question is
solved except for the two special cases of dimensions $n=4$ and $n=5$,
though in the case $n=5$ it is at least known that the statement is
true if only homeomorphisms are required.  (For the case $n=4$ this
question is equivalent to the Poincar\'e conjecture; see
\cite{Wils67}.)

\begin{prop} \label{Tconstruct}
  Let $n\neq 4$ and let $V:\R^n\to\R$ be a proper, positive definite
  $C^1$ function. Assume furthermore that $V$ is smooth on
  $\R^n\setminus\{0\}$ with nonvanishing gradient.  
  Then for each class $\KK_\infty$
  function $\gamma$ which is smooth on $(0,\infty)$ there exists
  a homeomorphism $T:\R^n \to \R^n$ with $T(0) = 0$ such that
\[ \tilde{V}(y) :=  V(T^{-1}(y)) = \gamma (\|y\|) \,.\]
In particular this holds for $\gamma(\|y\|) = \|y\|^2$.

If $n\neq 4,5$ then $T$ can be chosen to be a diffeomorphism on $\R^n
\setminus \{0\}$. Furthermore, in this case there exists a class
$\KK_\infty$ function $\gamma$ which is smooth on $(0,\infty)$ and
satisfies $\gamma(s)/\gamma'(s)\ge s$ such that $T$ is $C^1$ with
$DT(0)=0$.
\end{prop}
\begin{pf} 
  For the function $V$ the right hand side of the 
  {\em normed} gradient flow 
  \[ \dot{x} = \frac{\nabla V(x)'}{\|\nabla V(x)\|^2}\]
  is well defined and smooth for $x \ne 0$.  Denote the solutions by
  $\psi(t,x)$.  Then $V(\psi(t,x)) = V(x)+t$, and thus since $V$ is
  proper and $\nabla V(x)\ne 0$ for $x\ne 0$ for a given initial value
  $x\in\R^n$ $\psi$ is well defined for all $t\in(-V(x),\infty)$, thus
  also smooth (see e.g.\ \cite[Corollary 4.1]{Hart82}).

  Fix $c>0$. We define a map $\pi: \R^n \setminus \{0\} \to
  V^{-1}(c)$ by
  \[ \pi(x) = \psi(c-V(x),x) \,.\]
  Obviously $\pi$ is smooth, and since the gradient flow crosses each
  level set $V^{-1}(a), a>0$ exactly once it induces a diffeomorphism
  between each two level sets of $V$, which are $C^\infty$ manifolds
  due to the fact that $V$ is smooth away from the origin with
  nonvanishing gradient.
  
  Now observe that the properties of $V$ imply that $V^{-1}(c)$ is a
  homotopy sphere (cf.\ also \cite[Discussion after Theorem
  1.1]{Wils67}), which implies that $V^{-1}(c)$ is diffeomorphic to
  $\Sp^{n-1}$ for $n=1,2,3$ (see e.g.\ \cite[Appendix]{Miln65a} for
  $n=2$, \cite[Theorem 3.20]{DuFN90} for $n=3$; $n=1$ is trivial).
  For $n\ge 6$ we can use the fact that the sublevel set
  $\{x\in\R^n\,|\,V(x)\le c\}$ is a compact, connected smooth manifold
  with a simply connected boundary, which by \cite[\S{}9, Proposition
  A]{Miln65} implies that the sublevel set is diffeomorphic to the
  unit disc $D^n$, hence $V^{-1}(c)$ is diffeomorphic to $\Sp^{n-1}$.
  Thus for all dimensions $n\neq 4,5$ we may choose a diffeomorphism
  $S: V^{-1}(c) \to \Sp^{n-1}$. By \cite{Free82} we may choose $S$ to
  be at least a homeomorphism in the case $n=5$.
 
  Let $Q:=S\circ\pi$. The coordinate transformation $T$ is now given by
$T(0)
  = 0$ and
  \[ T(x) =  \gamma^{-1}(V(x))\,Q(x)\,,\quad x\neq 0\,.\]
  An easy computation verifies that $T^{-1}(0) = 0$ and 
  \[ T^{-1}(y) = \psi\left(\gamma(y)-c, 
    S^{-1}\left(\frac{y}{\|y\|}\right)\right)\,,\quad y\neq 0\,,\]
  hence $T$ is a diffeomorphism on $\R^n\setminus\{0\}$ (resp. a
  homeomorphism if $n=5$).  Since $V(0)=0$, and
  $\psi(t,S^{-1}({y}/{\|y\|}))\to 0$ as $t\searrow -c$, both $T$ and
  $T^{-1}$ are homeomorphisms.

  Finally, we have that
  \begin{eqnarray*} 
  V(T^{-1}(y)) & = & V\left(\psi\left(\gamma(\|y\|)-c, 
  S^{-1}\left(\frac{y}{\|y\|}\right)\right)\right) \\ & = & 
  V\left(S^{-1}\left(\frac{y}{\|y\|}\right)\right) - c + \gamma(\|y\|) =
  \gamma(\|y\|)  
  \end{eqnarray*}
  which finishes the proof of the first assertion. 

  For $n\neq 4,5$ and $s>0$ we define 
  \[L(s):=\sup_{V(x)=s}\| DQ(x)\|\]
  and choose any class $\KK$ function $a$ which is $C^\infty$ and satisfies
  \[ a(s) \le \frac{s}{L(s)} \frall
  s\in(0,1] \,.\] 
  Then the function $h$ given by
  \[
     h(r)  = \int_0^r a(s) ds
  \]
  is smooth and of class $\KK_\infty$. 
  Note that this construction implies $h(r)\le r a(r)$ for all $r \ge
0$,
  hence $h(r)/h'(r)\le r$. Thus $\gamma:=h^{-1}$ is of class
$\KK_\infty$,
  smooth on $(0,\infty)$, and satisfies
  \[ \frac{\gamma(s)}{\gamma'(s)} = h^{-1}(s)h'(h^{-1}(s)) \ge 
  \frac{h^{-1}(s)h(h^{-1}(s))}{h^{-1}(s)} = s .\]
     
  Differentiating $T$ yields
  \[ DT(x) = h'(V(x))Q(x)\cdot DV(x) + h(V(x))DQ(x). \]
  For $x\to 0$ the first term tends to 0 since both $h'(V(x))=a(V(x))$
and
  $DV(x)$ tend to 0, and the second tends to 0 since for all $x$
sufficiently
  close to 0 the inequality
  \[ h(V(x)) \|DQ(x) \|\le a(V(x)) \|DQ(x)\| \le
  \frac{V(x)}{L(V(x))} \| DQ(x)\| \le V(x) \] holds by construction of
$h$.
  Thus $DT(x)\to0$, as $x\to 0$, and consequently $T\in C^1$ with
$DT(0)=0$,
by a straightforward application of the mean value theorem, see e.g.
\cite[Chap. V, Theorem 3.2]{Lang86} and the fact that a function is
continuously differentiable if all partial derivatives exist and are
continuous.  \qed
  \end{pf}

\section{Main Results}
Using the coordinate transformation $T$ we can now prove our main
results. 

\begin{thm} \label{gas=ges}
  Let $n\ne 4,5$ and consider any system \eqref{system} on $\R^n$
  which is UGAS \eqref{gas}.  We suppose that the set $D\subset \R^m$
  is compact.  Then,
  \eqref{system} can be transformed into a system
  \eqref{transform} that is UGES \eqref{ges}.
  
  In particular, the constants in \eqref{ges} can be chosen to be $c=1,
  \lambda=1$.
\end{thm}
\begin{pf}
Under our assumptions, by \cite[Theorem 2.9, Remark 4.1]{LiSW96}\footnote{%
To be precise, the results in that reference make as a blanket assumption the
hypothesis that $f$ is locally Lipschitz, not merely continuous, at $x=0$.
However, as noted in e.g.~\cite{w96}, the Lipschitz condition at the
origin is not used in the proofs.}
there exists a smooth function 
$V : \R^n \to \R$ for \eqref{system} such that
\begin{equation}
  \label{vdot}
   L_{f_d}V(x) \le -\alpha_1(\|x\|)
\end{equation}
for some class $\KK_\infty$ function $\alpha_1$. 
Furthermore,
there exist class $\KK_\infty$ functions $\alpha_2,\alpha_3$ such
that
 \begin{equation}
   \label{bounds}
   \alpha_2(\|x\|) \le V(x) \le \alpha_3(\|x\|)\,.
 \end{equation}
 
 Now let $\alpha_4$ be a $C^1$ function of class $\KK_\infty$ which is
 smooth on $(0,\infty)$ and satisfies $\alpha_4'(0)=0$, such that
 $\alpha_4(a) \le \min\{a,\alpha_1\circ\alpha_3^{-1}(a)\}$ for all
 $a\ge0$.
 
 Such a function can be obtained e.g.\ by a slight modification of the
 construction in \cite[Proof of Lemma 11]{PraW96}: Take a class
 $\KK_\infty$ function satisfying 
 $\delta(a)
 \le
 \min\{a,\alpha_1\circ\alpha_3^{-1}(a)\}$ and which is smooth on
 $(0,\infty)$. Then
\[ \alpha_4(a) = \frac{2}{\pi}\int_0^a
\frac{\delta(\tau)}{1+\tau^2}d\tau\]
has the desired properties. Thus we obtain
\begin{equation}
  \label{vdot2}
   L_{f_d}V(x) \le -\alpha_4(V(x)).
\end{equation}
Now define 
\[ \rho(a) := \exp\left(-\int_a^1  
\alpha_4(\tau)^{-1} d\tau\right) \mbox{ for }a>0,\quad \rho(0):=0\]
Obviously $\rho$ is smooth on $(0,\infty)$; furthermore 
$\rho$ is of class $\KK_\infty$ and by \cite[Lemma 12]{PraW96} $\rho$ is 
a $C^1$ function on $[0,\infty)$ with $\rho'(0)=0$. 
Thus defining
\[ W(x) := \rho(V(x)) \]
we obtain a $C^1$ Lyapunov function, which is smooth on
$\R^n\setminus\{0\}$,
for which an easy calculation shows that
\[ L_{f_d}W(x) = \frac{\exp\left(-\int_{V(x)}^1  
    \alpha_4(\tau)^{-1} d\tau\right)}{\alpha_4(V(x))} L_{f_d}V(x) \le
-W(x).
\] 
Applying Proposition \ref{Tconstruct} to $W$, 
using the class $\KK_\infty$
function $\gamma$ with $\gamma(s)/\gamma'(s)\ge s$
we obtain for each $d\in
D$ and $y\neq 0$
\[ \langle \tilde{f}(y,d), y \rangle =
\frac{\|y\|}{\gamma'(\|y\|)}L_{\tilde{f}_d}\tilde{W}(y)  
\le -\frac{\|y\|}{\gamma'(\|y\|)}\tilde{W}(y) = 
-\frac{\|y\|}{\gamma'(\|y\|)}\gamma(\|y\|) 
\le - \|y\|^2\,. \]
Clearly the overall inequality also holds for $y=0$ so that we obtain
\[ \frac{d}{dt}\|y(t)\|^2 = 2 \langle \tilde{f}(y(t),d(t)), y(t) \rangle
\le - 2\|y(t)\|^2 \] and hence $\|y(t)\|^2 \le e^{-2t}\|y(0)\|^2$, i.e.\
the
desired exponential estimate. \qed
\end{pf}

\begin{thm} \label{iss=ises}
  Let $n\ne 4,5$ and suppose that the system \eqref{system} on $\R^n$
  is ISS \eqref{iss} with some class $\KK_\infty$ function $\alpha$
  and some class $\KK\LL$ function $\beta$.  Then \eqref{system} is
  can be transformed into a system \eqref{transform} that is ISES \eqref{ises}
with
  constants $c=\lambda=1$.
\end{thm}
\begin{pf}
By \cite[Theorem 1]{SonW95}\footnote{%
As with the UGAS proof, it is easy to verify that the assumption that the
right-hand side is Lipschitz at zero is never actually used in~\cite{SonW95}.
The possible non-uniqueness of trajectories does not affect the argument used
in Lemma 2.12 in that paper, which reduces the problem to one of UGAS.}
there exists a $C^1$ function $V$ 
which is smooth on $\R^n\setminus\{0\}$ 
and a class $\KK_\infty$ function $\chi$ such that 
\[ \|x\| > \chi(\|d\|) 
\quad \Rightarrow \quad L_{f_d}V(x) \le -\alpha_1(\|x\|) 
\] 
for some class $\KK_\infty$ function $\alpha_1$. 
Furthermore,
there exist class $\KK_\infty$ functions $\alpha_2,\alpha_3$ such
that
 \[
   \alpha_2(\|x\|) \le V(x) \le \alpha_3(\|x\|)\,.
 \]
As in the proof of Theorem \ref{gas=ges}
we find a 
function $\rho$ which is class $\KK_\infty$, $C^1$, and smooth on
$\R^n\setminus\{0\}$,
such that $W = \rho\circ V$ satisfies
\[ \|x\| > \chi(\|d\|) \quad \Rightarrow \quad L_{f_d}W(x) \le -W(x)\,.
\] 
Now Proposition \ref{Tconstruct} yields a parameter transformation
$T$ such that $\tilde{W}(y)=W(T^{-1}(y)) = \gamma(\|y\|)$ and 
$\gamma(s)/\gamma'(s)\ge s$. 

Now choose a class $\KK_\infty$ function $\delta$ such that
$\|T^{-1}(y)\|\ge \delta(\|y\|)$ and define $\tilde{\alpha} =
\delta^{-1}\circ \chi$. Then a straightforward calculation yields
\begin{equation} \|y\| > \tilde{\alpha}(\|d\|) \quad 
\Rightarrow \quad L_{\tilde{f}_d}\tilde{W}(y) \le -\tilde{W}(y).
\label{quadgain}\end{equation} 
Similar to the proof of Theorem \ref{gas=ges} 
this implies 
\[ \|\tilde{\phi}(t,y,d(\cdot))\| \le e^{-t}\|y\|\] as long as 
$ \|\tilde{\phi}(t,y,d(\cdot))\|> \tilde{\alpha}(\supdt)$
which yields the desired estimate. \qed
\end{pf}

\begin{thm} \label{ises=hinf}
  Consider the system \eqref{system} on $\R^n$ being ISES \eqref{ises}
  with some class $\KK_\infty$ function $\alpha$
  and $c=\lambda=1$.  
  Then there exists a
  homeomorphism $R:\R^m \to \R^m$ on the input space with $R(0) = 0$,
  that is a diffeomorphism on $\R^m \setminus \{0\}$, such that the
  the transformed system \eqref{dtransform} satisfies the nonlinear
  $H_\infty$ estimate \eqref{hinf}.
\end{thm}
\begin{pf}
{}From \eqref{ises} it is immediate that for any $d(\cdot)\in\DD$, 
any $x\in\R^n$, and any $T>0$ we have 
\begin{equation}
\|x\| \ge e^{{T}}\alpha(\supdT)
\quad 
\Rightarrow \quad  \|\phi(t,x,d(\cdot))\| \le e^{-{t}}\|x\|
\frall t\in[0,T] \,.
\label{qgain}\end{equation}
Now consider the function $W(x)=\|x\|^2$. Then \eqref{qgain} implies 
\[
\|x\| \ge e^{{T}}\alpha(\supdT)
\quad \Rightarrow \quad  W(\phi(t,x,d(\cdot))) \le e^{-2t}W(x)
\frall t\in[0,T]\,.
\]
In particular this estimate is valid for constant functions 
$d(\cdot)\equiv d\in D$,
thus the mean value theorem (observe $W(\phi(0,x,d)) = W(x) =
e^{-0}W(x)$)
yields
\[ 
\|x\| \ge \alpha(\|d\|) \quad 
\Rightarrow \quad  L_{{f}_d}{W}(x) \le -2W(x) \le -W(x)\,.
\] 
Now defining 
\[ \tilde{\alpha}(r) = \sup_{\|x\|\le \alpha(r),\|d\|\leq r} 
\langle f(x,d),x\rangle \]
we obtain a class $\KK_\infty$ function $\tilde{\alpha}$ with
\[ L_{f_d}{W}(x) \le -{W}(x) + \tilde{\alpha}(\|d\|). \]
Without loss of generality (one could take a larger $\tilde{\alpha}$),
we may assume $\tilde{\alpha}$ to be smooth on $(0,\infty)$,
and thus
\[ R(d) := \frac{\tilde{\alpha}(\|d\|)^2 d}{\|d\|}\]
has the regularity properties as stated in the assertion. Now the 
transformation \eqref{dtransform} yields 
\[  L_{\bar{f}_v}{W}(x) \le -{W}(x) + \|v\|^2.\]
Integrating this equation along a trajectory $x(\cdot)$ gives 
\[  W(x(t)) - W(x(0)) \le - \int_0^t W(x(s))ds + \int_0^t\|v(s)\|^2ds 
\]
which implies \eqref{hinf} since $W(x) = \|x\|^2$. \qed
\end{pf}

\section{Remarks}
Note that, in general, for our results to be true we cannot expect $T$
to be diffeomorphic on the whole $\R^n$.  Consider the simplest case
where $f$ does not depend on $d$ and is differentiable at the origin.
If $T$ were a diffeomorphism globally, then $DT^{-1}(0)$ would be
 well-defined, which implies that
\[ D\tilde{f}(0) = 
\left.\frac{\partial}{\partial
    y}\right|_{y=0}DT(T^{-1}(y))f(T^{-1}(y)) = DT(0) Df(0)
DT^{-1}(0)\] and so the linearizations in $0$ are similar; in
particular, the dimension of center manifolds remains unchanged.

Actually, if one wants the exponential decay to be $e^{-t}$, even for
linear systems one cannot obtain a diffeomorphism $T$.  As an example,
consider the one-dimensional system $\dot x = -x/2$.  Here one uses
the change of variables $y=T(x)$ given by $T(x)=x^2 ,\, x > 0,\ 
T(0)=0$ and $T(x)=-x^2 ,\, x < 0$ to obtain $\dot y = -y$. Note that
$T$ is $C^1$ with $DT(0)=0$. The inverse of this $T$ is given by
$T^{-1}(y)= \sqrt{y},\, y > 0,\ T^{-1}(0)=0$ and $T^{-1}(y)=-\sqrt{-y}
,\, y < 0$ which is smooth only away from the origin, though
continuous globally.

An example for the case of nontrivial center manifolds is given by the
system $\dot{x} = -x^3$. Let us first note that for this system there
is no transformation in the class we consider such that the
transformed system is of the form $\dot{y}=-y$. The reason for this is
that we would have $\dot{T}(x) = \dot{y} = -y = -T(x)$, so at least
for $x>0$ $V=T$ is a Lyapunov function with the property that
$\dot{V}(x) = -V'(x)x^3 = -V(x)$. It is readily seen that the
solutions of this differential equation (in $x$ and $V$) are $V_c(x) =
c \exp{{1\over - 2x^2}}$, for $c\in\R$. However, the image of
$[0,\infty)$ under such $V_c$ yields a bounded set, so that these
functions are no candidate for coordinate transforms on $\R$.
Nonetheless a coordinate transform according to our requirements can
now be easily built: Take any $\KK_\infty$ function $\alpha$ with
$\alpha'>0$ on $(0,\infty)$ so that with via the symmetrization
$\alpha(-x) := \alpha(x)$ we get a smooth function on $\R$. Now define
\[ T(x) := \alpha(x)V_1(x),\, x\ge 0,\quad 
T(x) := - \alpha(x)V_1(x),\, x <0\,.\] Then for $y\neq 0$ we have
$\dot{y} = \dot{T(x)} = -(1+ {\alpha'(x)x^3 \over \alpha(x)}) T(x) <
-y$, so that the transformed system decays at least exponentially with
constants $c=1,\lambda=1$. Again note that the requirement $DT(0)=0$
is vital, in fact all orders of derivatives vanish in $0$.

A basic ingredient of the proof of Theorem~\ref{gas=ges} is the
construction of a Lyapunov function with the property $\dot{V} \leq
-V$. Actually, one may even, under restricted conditions, obtain the
equality $\dot{V} = -V$. It should be noted that already in
\cite{Bhat67} it is shown that for dynamical systems with globally
asymptotically stable fixed point a {\em continuous} Lyapunov function
with the property $V(\phi(t,x)) = e^{-t}V(x)$ exists, see also Chapter
V.2 in \cite{BhatSzeg70}. Note, however that in these references only
systems with trajectories defined on $\R$ are considered, which does
not include the previous example.  Indeed, if $f(x,d)=f(x)$ is
independent of $d\in D$ and the system $\dot{x}=f(x)$ is backward
complete we can can also define a coordinate transformation based on a
different $W$ than the one used in the proof of Theorem~\ref{gas=ges}:
In this case the function $W(x)=\exp t(x)$ with $t(x)$ defined by
$V(\phi(t(x),x)) = 1$ is positive definite, proper, and satisfies
$L_{f}W(x) = -W(x)$, thus $W(\phi(t,x)) = W(x) - t$.  Since
$V^{-1}(1)=W^{-1}(1)$ we still find a diffeomorphism $S$ as in the
proof of Proposition \ref{Tconstruct}.  Deviating from this proof,
instead of the gradient flow we now use the trajectories of the
system, i.e.\ we define $\pi(x) = \phi(W(x)-1,x)$ yielding $W(\pi(x))
= W(x)-(W(x)-1) = 1$.  Thus from $\pi$ we can construct $T$ as in the
proof of Proposition \ref{Tconstruct}, and obtain $W(T^{-1}(y)) =
\|y\|^2$. Furthermore the definition of $\pi$ implies that each
trajectory $\{\phi(t,x)\,|\,t\in\R\}$ is mapped onto the line $\{
\alpha S(\pi(x))\,|\, \alpha>0\}$ and consequently $\tilde{f}(y)=-y$,
i.e.\ we obtain a transformation into the linear system $\dot{y}=-y$.
Note, however, that with this construction the coordinate
transformation will in general only have the regularity of $f$ (e.g.\ 
a homeomorphism if $f$ is only $C^0$), which is inevitable since it
transforms $f$ into a smooth map. Moreover, this construction cannot
be generalized to systems with disturbances.

Since we are not requiring that the inverse of a change of variables
be itself a change of variables (because one may, and in fact does in
our constructions, have $DT(0)=0$, in which case $T^{-1}$ is not
differentiable at the origin), the way to define a notion of
``equivalence'' is by taking the transitive and symmetric closure of
the relation given by such changes of variables.  That is, we could
say that system \eqref{system} is {\em equivalent\/} to a
system~\eqref{transform} if there exist $k\in\N$ and maps $f_0 = f,
f_1, \ldots, f_k =\tilde{f}:\R^n\times D \to \R^n$, all satisfying the
assumptions on $f$, with the following properties: For each
$i=0,\ldots,k-1$ there exists a change of variables $T$ as above such
that $f_l(y,d) = DT(T^{-1}(y)) f_m(T^{-1}(y),d)$, where $l=i, m=i+1$
or $l=i+1,m=i$.

Finally, regarding our notion of system transformation, note that even
if $f(0,d) \neq0$ for some $d\in\DD$ for the original system
(\ref{system}), then under the assumption $DT(0)=0$ we have
$\tilde{f}(0,d) = 0$ for all $d\in\DD$ for the transformed system.
This implies that even if the original system had unique trajectories
through zero, the transformed system cannot have this property.

{\bf Acknowledgments:} We thank David Angeli for suggestions regarding
the remark on backward complete systems, as well as Uwe Helmke and
Laurent Praly for many references to the literature.


\end{document}